
\input amstex
\documentstyle{amsppt}
\nologo
\magnification=\magstep1
\hsize=145mm
\vsize=220mm
\hcorrection{-7mm}
\vcorrection{-10mm}
\baselineskip=18pt
 \abovedisplayskip=4pt
 \belowdisplayskip=4pt
\parskip=3pt
\parindent=8mm


\define\dnl{\newline\newline}
\define\lra{{\longrightarrow}}
\define\nl{\newline}

\define\Bs{{\roman{Bs}}}
\define\Exc{{\roman{Exc}}}
\define\Frac{{\roman{Frac}}}
\define\Int{{\roman{Int}}}
\define\Max{{\roman{Max}}}
\define\Sup{{\roman{Sup}}}
\define\Supp{{\roman{Supp}}}
\define\dm{{\roman{dim}}}

\define\SO{\Cal O}
\define\BC{{\Bbb C}}
\define\BP{{\Bbb P}}
\define\BQ{{\Bbb Q}}

\define\Cor{{\bf Corollary. }}

\define\Exa{{\it Example. }}
\define\Fact{{\it Fact. }}
\define\Lem{{\bf Lemma. }}
\define\Prf{{\it Proof. }}

\define\Rmk{{\it Remark. }}
\define\Th{{\bf Theorem. }}

\topmatter
\author Takao FUJITA \endauthor
\address {Takao FUJITA\newline
Department of Mathematics
\newline
Tokyo Institute of Technology
\newline
Oh-okayama, Meguro, Tokyo
\newline
152 Japan
\newline
e-mail:fujita\@math.titech.ac.jp}
\endaddress

\title Towards a separation theorem of points \\
 by adjoint linear systems on polarized threefolds
\endtitle
\endtopmatter

\document

\noindent {\bf \S0. Introduction}

The motivation of this note comes from the following

(0.1) {\bf Conjecture A. } {\sl
Let $L$ be an ample line bundle on a smooth algebraic variety $M$.
Let $K$ be the canonical bundle of $M$ and set $n=\dm M$.
Then $K+tL$ is spanned for every $t>n$, i.e.,
$\Bs{\vert}K+tL{\vert}=\emptyset$.
Moreover, $K+nL$ is spanned unless $L^n=1$.}

(0.2) {\bf Conjecture B. } {\sl
Let $L$, $M$, $K$ be as above.
Then $K+tL$ is very ample for every $t\ge n+2$.
Moreover $K+(n+1)L$ is very ample unless $L^n=1$.}

(0.3) We have now many partial answers to the above conjectures,
but they depend essentially on vanishing theorems of Kodaira type.
Hence we consider the complex cases only, though no counter-example is known in
positive characteristic cases.

Now we review such partial results.

(0.4) \Fact {\sl
The conjectures A and B are true when $n=2$.}

This is indeed a consequence of the famous result of Reider [{\bf R}].

(0.5) \Fact {\sl
Conjecture A is true if $L$ is spanned.}

{\it Outline of proof.}
Let $\pi:M_1\lra M$ be the blow-up at $x\in M$ and let $E$ be the exceptional
divisor over $x$.
Then $\Bs{\vert}\pi^*L-E{\vert}$ has only finitely many points off $E$, since
${\vert}L{\vert}$ defines a finite morphism $\rho:M\lra \BP^N$.
Hence $\pi^*L-E$ is nef, so $H^1(M_1,\pi^*(K+nL)-E)=0$ by Kawamata-Viehweg
vanishing theorem.
This implies $x\notin\Bs{\vert}K+nL{\vert}$.
The cases $K+tL$ with $t>n$ is easier.

(0.6) \Fact {\sl
Conjecture B is true if $L$ is very ample.}

{\it Outline of proof.}
Let $\pi$, $M_1$, $E$ be as above.
It suffices to show $\Bs{\vert}\pi^*(K+tL)-E{\vert}=\emptyset$.
If not, take a base point $y$ of ${\vert}\pi^*(K+tL)-E{\vert}$.
Then $\pi(y)$ is a point on $M$ which may be infinitely near to $x$.
Let $\ell$ be the line passing $\pi(y)$ and $x$ in $\BP$ with respect to the
embedding $\rho: M\lra\BP$ defined by ${\vert}L{\vert}$.

If $\ell\not\subset M$, let $\pi_2: M_2\lra M_1$ be the blow-up at $y$ and let
$E_2$ be its exceptional divisor.
Then $L_{M_2}-\pi_2^*E_1-E_2$ is nef as in (0.5), and we get a contradiction by
using vanishing theorem.

If $\ell\subset M$, there is a smooth member $D$ of ${\vert}L{\vert}$ such that
$\ell\subset D$, provided $n\ge 3$.
Let $D_1$ be the proper transform of $D$ on $M_1$.
Then $y\in D_1$ and $[\pi^*(K+tL)-E_1]_{D_1}=[\pi^*(K_D+(t-1)L_D)-E_1]_{D_1}$
is spanned at $y$ by the induction hypothesis on $n$.
This yields a contradiction, as desired.

(0.7) \Fact{\sl
$K+tL$ is nef under the hypothesis of Conjecture A.}
(cf., e.g., [{\bf F2}]).

(0.8) \Exa
Let $(M,L)$ be a weighted hypersurface of degree 6 in the weighted projective
space $\BP(3,2,1,\cdots,1)$.
Then $K=(1-n)L$, $L^n=1$, $\Bs{\vert}L{\vert}\ne\emptyset$, $2L$ is spanned but
not very ample, and $mL$ is very ample for every $m\ge3$ (cf., e.g., [{\bf
F0};(6.14-17)]).
This example shows that the bounds for $t$ in (0.1) and (0.2) are sharp.

(0.9) \Fact{\sl
Conjecture A is true for $n=3$.}

This is a consequence of the remarkable result of Ein-Lazarsfeld [{\bf EL}],
supplemented by [{\bf F4}].

(0.10) Now we focus our attention to Conjecture B in case $n=3$.
In this note we describe an approach to this problem, based on a result in
[{\bf F3}].
This method is less powerful than Ein-Lazarsfeld's, but it is cheaper, i.e.,
needs fewer computations and pages.
The main result can be stated as follows.

(0.11) \Th{\sl
Let $(M,L)$ be a smooth polarized threefold over $\BC$.
Then ${\vert}K+tL{\vert}$ separates any two different points on $M$ for
$t\ge6$.}

The bound for $t$ is not the one predicted in (0.2), and at present I cannot
prove the separation of infinitely near points.
But I hope that the method is of some interest because of its cheapness.
I remark also that the same method gives the spannedness of
${\vert}K+tL{\vert}$ for $t\ge5$.
\dnl
{\bf \S1. Preliminaries}

(1.1) {\it Notation.}
Given a $\BQ$-divisor $D$, we denote the integral part and fractional part of
$D$ by $\Int(D)$ and $\Frac(D)$ respectively.
Thus $D=\Int(D)+\Frac(D)$ and $\Frac(D)$ is effective.
$D$ is called an {\it NC $\BQ$-divisor} if its support has only normal crossing
singularities.

Given a birational morphism $f:V\lra W$, we denote by $\Exc(f)$ (or
$\Exc(V/W)$) the exceptional set of $f$.
Note that $\Exc(f)=f^{-1}(X)$ for the smallest subset $X$ of $W$ such that
$f^{-1}(W-X)\cong W-X$.

(1.2) \Lem{\sl
Let $\pi:M_1\lra M$ be the blow-up at a smooth point $x$ on a variety $M$ with
$\dm M=n$ and let $E$ be the exceptional divisor.
Then, for any line bundle $L$ on $M$,
$h^0(\pi^*L-tE)\ge h^0(L)-\frac1{n!}t(t+1)\cdots(t+n-1)$.}

\Prf
By the exact sequence
$0\lra H^0(\pi^*L-(s+1)E)\lra H^0(\pi^*L-sE)\lra H^0(E,-sE_E)$
we infer
$h^0(\pi^*L-sE)-h^0(\pi^*L-(s+1)E)\le h^0(\BP^{n-1},\SO(s))$.
Therefore
$h^0(L)-h^0(\pi^*L-tE)\le \sum_{s=0}^{t-1}\frac1{(n-1)!}(s+1)\cdots(s+n-1)
=\frac1{n!}t(t+1)\cdots(t+n-1)$.

(1.3) Now we recall the following

\Th {\sl
Let $L$ be a line bundle on a variety $V$ with $\dm V=n$.
Suppose that $h^0(V,tL)\ge dt^n/n!+\phi(t)$ for infinitely many positive
integers $t$, where $d$ is a positive number and $\phi(t)$ is a function such
that
$\lim_{t\rightarrow\infty}\phi(t)/t^n=0$.
Then, for any $\epsilon>0$, there is a birational morphism $\pi: M\lra V$
together with an effective $\BQ$-divisor $E$ on $M$ such that $H=\pi^*L-E$ is a
semiample $\BQ$-bundle with $H^n>d-\epsilon$.}

For a proof, see [{\bf F3}].

(1.4) \Rmk
In the above situation, we can choose $M$ and $E$ so that\nl
1) $H$ is ample,\nl
2) $E$ is an NC $\BQ$-divisor, and\nl
3) $\Exc(\pi)$ is contained in the support of $E$.

Indeed, we replace $M$ and $E$ as follows.
Since $H$ is nef and big, there is an effective $\BQ$-divisor $\Delta$ such
that $H-\delta\Delta$ is ample for any $0<\delta\le 1$.
Let $f: M'\lra M$ be a succession of blowing-ups along smooth centers such that
${\roman{Supp}}(f^*(E+\Delta))\cup\Exc(\pi\circ f)$ is an NC divisor.
Let $N$ be an effective $\BQ$-divisor on $M'$ such that $\Supp(N)=\Exc(f)$ and
$-N$ is $f$-ample.
Then $f^*(H-\delta\Delta)-\alpha N$ is ample for $0<\alpha\ll 1$.
Set $E'=f^*(E+\delta\Delta)+\alpha N$
for sufficiently small positive numbers $\delta$ and $\alpha$.
Then this pair $(M', E')$ satisfies the conditions 1) and 2), preserving the
property $(H')^n>d-\epsilon$.
Adding a very small portion of $\Exc(\pi\circ f)$ to $E'$, we can arrange
things so that 3) is satisfied too, since the requirements are open conditions.

(1.5) {\bf Index Theorem. }{\sl
Let $A$, $B$, $H$ be nef $\BQ$-bundles on a variety $V$ with $\dm V=n$.
Then $(ABH^{n-2})^2\ge(A^2H^{n-2})(B^2H^{n-2})$.}

For a proof, reduce the problem to the case in which $H$ is very ample, and
then use the induction on $n$.
See, e.g., [{\bf F1};(1.2:4)] for details.

(1.6) \Cor {\sl
Let $A$, $B$ be nef $\BQ$-bundles on a threefold $V$.
Then $(A^3)^2(B^3)\le(A^2B)^3$.}
\dnl
{\bf \S2. Main Result}

(2.1) Throughout this section, let $A$ be a line bundle on a threefold $V$
having only log terminal singularities and $K$ be the canonical $\BQ$-bundle of
$V$.
Let $x$, $y$ be different smooth points on $V$.
Suppose that $B=A-K$ is nef and satisfies the following conditions.\nl
1) $BC\ge d_1>2$ for any irreducible curve $C$ passing $x$ or $y$,\nl
2) $B^2S\ge d_2$ for any irreducible surface $S$ containing $x$ or $y$,\nl
3) $B^3\ge d_3$,\nl
where $d_i$'s are positive numbers satisfying\nl
NA1) $d_3^2(d_3(d_1-1)^3-54d_1^3)>(d_3(d_1-1)-d_1d_2)^3>0$.

Note that these conditions are satisfied by $d_i=6^i$ when $B=6L$ for some
ample line bundle $L$.

{\it Remark } A.
We may assume that $d_i\in\BQ$ by perturbing slightly if necessary.

{\it Remark } B.
We may assume $d_3=B^3$.

To see this, set $q=(d_1-1)/d_1$ and $b=1-(d_2/qd_3)$.
Then NA1) is equivalent to $d_3^2(q^3d_3-54)>(qd_3-d_2)^3>0$, and further to
$54<q^3d_3-d_3(q-d_2/d_3)^3=q^3d_3(1-b^3)=q^2d_2(1+b+b^2)$.
Replacing $d_3$ by $d=B^3$, $b$ becomes larger and NA1) is still satisfied.

{}From now on, we assume $d_3=B^3=d$.

(2.2) Let $\pi_2: V_2\lra V$ be the blow-up at $x$ and $y$, and let $E_x$,
$E_y$ be the exceptional divisors over $x$, $y$ respectively.
Take a positive integer $m$ such that $mqB$ is Cartier, where $q=(d_1-1)/d_1$.
Then $h^0(tmqB)=\frac{t^3}6 m^3q^3d+O(t^2)$ for $t\gg0$.
Hence $h^0(tmqB_{V_2}-3tmE_x-3tmE_y)\ge\frac{t^3}6m^3(q^3d-54)+O(t^2)$ by
(1.2).
Since $q^3d-54>d(q-d_2/d)^3$ by NA1), there is a birational morphism $f: M\lra
V_2$ and an effective NC $\BQ$-divisor $F$ on $M$ such that
$H=(qB-3E_x-3E_y)_M-F$ is ample, $H^3>d(q-d_2/d)^3$ and $\Exc(\pi_2\circ
f)\subset\Supp(F)$, by virtue of (1.3) and (1.4).

(2.3) {\it Claim}. {\sl
Let $F+f^*(3E_x+3E_y)=\sum_i\mu_iF_i$ be the prime decomposition.
Then $\mu_i<1$ if $S=\pi(F_i)$ is a surface passing $x$ or $y$, where
$\pi=\pi_2\circ f$.}

\Prf
If not, we have $B_M^2(qB_M-H)\ge B_M^2F\ge B^2S\ge d_2$, so
$B_M^2H\le qd-d_2$.
On the other hand $d^2H^3\le (B_M^2H)^3$ by (1.6), hence
$(qd-d_2)^3\ge d^2H^3$, contradicting the choice of $H$.

(2.4) Let $R=\sum_ia_iF_i$ be the ramification divisor of $\pi=\pi_2\circ f$.
Set $\nu_i=\mu_i/(a_i+1)$ for each $i$.
Perturbing the choice of $F$ if necessary, we may assume that $\nu_i\ne\nu_j$
if $i\ne j$.
Let $I_x=\{i{\vert}\pi(F_i)\ni x\}$, $\nu_x=\Max_{i\in I_x}\nu_i$ and let $i_x$
be the index $i\in I_x$ such that $\nu_i=\nu_x$.
Replacing $x$ by $y$, we define $I_y$, $\nu_y$ and $i_y$ similarly.
Since $\mu_i\ge3$ and $a_i=2$ for the proper transform $F_i$ of $E_x$ or $E_y$,
we have $\nu_x\ge 1$ and $\nu_y\ge 1$, so $F_{i_x}$ and $F_{i_y}$ are
$\pi$-exceptional by (2.3).

(2.5) It is possible that $i_x=i_y$.
In this case $F_{i_x}=F_{i_y}$ is mapped onto a curve passing $x$ and $y$.
This case will be treated in (2.6).
Here we consider the case $i_x\ne i_y$.

By symmetry we may assume $\nu_x>\nu_y$.
Put $c=\nu_y^{-1}$ and $W=F_{i_y}$.
Then $c\mu_i\le a_i+1$ for any $i\in I_y$, with the equality only for $i=i_y$.
Hence $\Int(\sum_i(c\mu_i-a_i)F_i)=W+P-N$ for some effective divisors $P$, $N$
with $y\notin\pi(P)$.
Note that $F_{i_x}$ is a component of $P$ since $\nu_x>\nu_y$.
Set $\Delta=\Frac(\sum_i(c\mu_i-a_i)F_i)$.
Then $A-K(M)-\Delta=B-R-\Delta=\frac1{d_1}B+(1-c)qB+cH+P-N+W$.
Hence $H^1(M,A-P+N-W)=0$ by the vanishing theorem, so
$H^0(M,A-P+N)\lra H^0(W,(A-P+N)_W)$ is surjective.

For any component $F_i$ of $N$, we have $a_i>c\mu_i>0$.
Hence $N$ is $\pi$-exceptional, so $H^0(M,A-P)\cong H^0(M,A-P+N)$.
Since $H^0(W,A-P)\lra H^0(W,A-P+N)$ is injective, we infer that
$H^0(M,A-P)\lra H^0(W,A-P)$ is surjective.

Now we divide the cases according to $\dm\pi(W)$.

(2.5.0) When $\dm\pi(W)=0$, we have $W\subset\pi^{-1}(y)$ and $P_W=0$.
Therefore $H^0(W,A-P)\cong\BC$.
Take $\varphi\in H^0(M,A-P)$ which is mapped to $1\in H^0(W,A-P)$.
This corresponds to $\psi\in H^0(M,A)$ vanishing along $P$, and $\psi=\pi^*\xi$
for some $\xi\in H^0(V,A)$.
The above choice of $\varphi$ implies $\xi(y)\ne 0$, while $\xi(x)=0$ since
$\pi(P)\supset\pi(F_{i_x})\ni x$.
Thus $x$ and $y$ are separated by ${\vert}A{\vert}$.

(2.5.1) When $C=\pi(W)$ is a curve, let $W\lra Z\lra C$ be the Stein
factorization.
Then $Z$ is normal and every fiber of $g: W\lra Z$ is connected.
Take a point $z$ on $Z$ lying over $y$ and let $\Gamma$ be the fiber
$g^{-1}(z)$.
Since $BC\ge d_1$, $d_1^{-1}B-\Gamma$ is nef on $W$.
Hence $(A-P+N-\Gamma)_W-K(W)-\Delta=(d_1^{-1}B-\Gamma)+(1-c)qB+cH$ is ample on
$W$, so
$H^1(W,A-P+N-\Gamma)=0$ and $H^0(W,A-P+N)\lra H^0(\Gamma,A-P+N)$ is surjective.
Argueing similarly as in [{\bf EL}], we infer that $\varphi_\Gamma\ne 0$ for
some $\varphi\in H^0(M,A-P)$.
This implies $\xi(y)\ne 0$ while $\xi(x)=0$ as in (2.5.0).

(2.6) Now we consider the case $i_x=i_y$.
We need the following assumption.

NA2) $d_3^2(d_3(d_1-2)^3-27d_1^3)>(d_3(d_1-2)-d_1d_2)^3>0$.
\nl
As NA1), this is satisfied when $d_i=6^i$.
We may assume $d_i\in\BQ$ and $d_3=B^3=d$ as before.
We set $r=(d_1-2)/d_1$.

(2.6.1) Let $\pi_1: V_1\lra V$ be the blow-up at $x$ and let $E_x$ be the
exceptional divisor.
As in (2.2), we have
$h^0(tmrB-3tmE_x)\ge\frac{t^3}6m^3(r^3d-27)+O(t^2)$ by (1.2).
Hence there is a birational morphism $f_1:M'\lra V_1$ and an effective NC
$\BQ$-divisor $F'$ on $M'$ such that $H'=rB-3E_x-F'$ is ample on $M'$ and
$(H')^3>d(r-d_2/d)^3$ by NA2).

(2.6.2) By replacing the model $M$ and $M'$ suitably, we may assume that $M'=M$
and $F+F'$ is an NC $\BQ$-divisor.

Indeed, we have a common model $M''$ such that $F_{M''}+F'_{M''}+\Exc(M''/V)$
is an NC $\BQ$-divisor.
Perturbing similarly as in (1.4), we replace $F_{M''}$ and $F'_{M''}$ by $F_2$
and $F_1$ so that all the requirements are satisfied.
In particular $H_2=qB-3E_x-3E_y-F_2$ and $H_1=rB-3E_x-F_1$ are ample on $M''$.

{}From now on we assume $M=M'=M''$ for the sake of brevity.

(2.6.3) Let $F_1+f_1^*3E_x=\sum_i\mu'_iF_i$ and set
$\nu'_i=\mu'_i/(a_i+1)$,
where $a_i$ is as in (2.4).
$F_i$ is non-exceptional iff $a_i=0$, and for such $i$ with $\pi(F_i)\ni x
\text{ or }y$ we have $\nu'_i=\mu'_i<1$ by NA2) as in (2.3).

(2.6.4) Let $I_x$ be as in (2.4), set $\nu'_x=\Max_{i\in I_x}\nu'_i$ and let
$i'_x$ be the index $i\in I_x$ such that $\nu'_i=\nu'_x$.
Then $\nu'_x\ge 1$ as in (2.4).
We define similarly $I_y$, $\nu'_y$ and $i'_y$.

(2.6.5) When $\nu'_y\ge 1$, we divide the cases as before, and show that
${\vert}A{\vert}$ separates $x$ and $y$ in case $i'_x\ne i'_y$.

In the case $i'_x=i'_y$, set $W=F_{i'_x}$ for brevity and put
$c'=(\nu'_x)^{-1}$.
Then $\Int(\sum_i(c'\mu'_i-a_i)F_i)=W+P'-N'$ for some effective divisors $P'$,
$N'$ with $\{x,y\}\notin\pi(P')$.
Similarly in (2.5) we infer that
$H^0(M,A-P'+N')\lra H^0(W,A-P'+N')$ and
$H^0(M,A-P')\lra H^0(W,A-P')$ are surjective.
Let $Z$ be as in (2.5.1) and take points $z_1$ and $z_2$ lying over $x$ and $y$
respectively.
Let $\Gamma_j$ be the fiber over $z_j$.
Then $\frac2{d_1}B-\Gamma_1-\Gamma_2$ is nef on $W$ and
$A-P'+N'-\Gamma_1-\Gamma_2-K(W)-\Delta'=\frac2{d_1}B-\Gamma_1-\Gamma_2+(1-c')rB+c'H_1$ is ample on $W$, so
$H^0(W,A-P'+N')\lra H^0(\Gamma_1\cup\Gamma_2,A-P'+N')$ is surjective.
{}From this we infer that ${\vert}A{\vert}$ separates $x$ and $y$ in this case.

(2.6.6) It remains the case $\nu'_y<1$.
In this case, for $0\le s\le 1$, we set
$\mu_i(s)=(1-s)\mu_i+s\mu'_i$ and
$\nu_i(s)=(1-s)\nu_i+s\nu'_i=\mu_i(s)/(a_i+1)$.
Set $\nu_x(s)=\Max_{i\in I_x}\nu_i(s)$ and
$\nu_y(s)=\Max_{i\in I_y}\nu_i(s)$.
Then $\nu_x(0)=\nu_y(0)>1$, $\nu_y(1)<1$ and $\nu_x(s)\ge 1$ for any $s$.

Let $s_0=\Sup\{s{\vert}\nu_x(s)=\nu_y(s)\}$.
Then we have $i_1, i_2\in I_x$ such that $i_1\notin I_y$, $i_2\in I_y$,
$\nu_{i_1}(s_0)=\nu_x(s_0)=\nu_y(s_0)=\nu_{i_2}(s_0)$ and
$\nu'_{i_1}>\nu'_{i_2}$.
The index $i_2$ may be different from $i_x=i_y$, but we set here $W=F_{i_2}$
for simplicity.
Put $c=\nu_x(s_0)^{-1}$.
Then $\Int(\sum_i(c\mu_i(s_0)-a_i)F_i)=W+P_0-N_0$ for some effective divisors
$P_0$, $N_0$ with $y\notin\pi(P_0)$.
But $F_{i_1}$ is a component of $P_0$.
Let $\Delta_0$ be the fractional part.
Then
$A-K(M)-\Delta_0=\frac{1+s}{d_1}B+(1-c)\frac{d_1-1-s}{d_1}B+c((1-s)H_2+sH_1)+W+P_0-N_0$, so
$H^0(M,A-P_0+N_0)\lra H^0(W,A-P_0+N_0)$ and
$H^0(M,A-P_0)\lra H^0(W,A-P_0)$ are surjective.
Moreover, as in (2.5.1), we have $\varphi\in H^0(M,A-P_0)$ such that
$\varphi(y)\ne 0$.
Hence $\xi(y)\ne0$ for the section $\xi\in H^0(V,A)$ corresponding to
$\varphi$, but $\xi(x)=0$ since $F_{i_1}$ is a component of $P_0$ and $i_1\in
I_x$.
Thus $x$ and $y$ are separated by ${\vert}A{\vert}$.

(2.7) From the above arguments we obtain the following

\Th {\sl Let $A$ be a line bundle on a threefold $V$ having only log terminal
singularities with canonical $\BQ$-bundle $K$.
Let $x$, $y$ be two smooth points on $V$.
Suppose that $B=A-K$ is nef and satisfies the following conditions:\nl
1) $BC\ge d_1>2$ for any irreducible curve $C$ passing $x$ or $y$,\nl
2) $B^2S\ge d_2$ for any irreducible surface $S$ containing $x$ or $y$,\nl
3) $B^3\ge d_3$,\nl
where $d_i$'s are positive numbers satisfying the numerical conditions NA1) and
NA2).
Then $x$ and $y$ are separated by ${\vert}A{\vert}$.}

(2.8) \Rmk By the same method, we can prove the following result.

Let $(V,D)$ be a log variety having only log terminal singularities.
Let $x$, $y$ be two smooth points on $V$, let $A$ be a line bundle on $V$ such
that $B=A-K(V,D)$ satisfies the same conditions as above.
Then ${\vert}A{\vert}$ separates $x$ and $y$.

(2.9) Any way, applying (2.7) for $d_i=6^i$, we get the following

\Cor {\sl Let $L$ be an ample line bundle on a smooth threefold $M$ with
canonical bundle $K$.
Then ${\vert}K+6L{\vert}$ separates any two different points on $M$.}

\enddocument